\documentclass[12pt]{amsart}
\usepackage{amsmath,amsfonts,amssymb,amscd,amsthm,amsbsy,epsf}
\newtheorem{thm}{Theorem}[section]

\newtheorem{cor}[thm]{Corollary}
\newtheorem{prop}[thm]{Proposition}

\def\D{{\mathcal D}}
\def\Re{\operatorname{Re}}
\def\Vol{\operatorname{vol}}
\def\iint{\int\!\!\!\int}
\begin{document}
\title{Moduli space theory for constant mean curvature surfaces immersed 
in space-forms}
\author{Alexandre Gon\c{C}alves and Karen Uhlenbeck}

\maketitle

\section{Introduction}

The study of constant mean curvature surfaces in a space-form has been an
active field since the work of H. Hopf  in the 1920's  and H.Liebmann in
the years around 1900. The questions which are generally of interest are
global questions of existence and uniqueness in complete 3-manifolds.
 We deal in this short
paper on a question of existence and uniqueness with respect
to the complex structure and the quadratic Hopf differential of a compact
surface
in a constant curvature 3-manifold which is not necessarily complete.Our
final result applies only to the case of surfaces
embedded
in a local 
3-dimensional space of constant curvature $-1$, where the mean curvature
constant $c$ satisfies $|c|< 1$. 
However, the technique suggests some approaches to the more interesting
cases, for example the work of Bryant \cite{B} and Kenmotsu \cite{K}.

The Gauss-Codazzi equations for constant mean curvature immersions of a
surface into a 3-dimensional space-form are a $3\times 3$ system of partial
differential equations of mixed order. Once a complex structure is
chosen, the equations break down into two equations. The Codazzi equation
on  the second fundamental 
form  yields the Cauchy-Riemann equation for a holomorphic
quadratic differential first noticed  and used by Heintz Hopf \cite{H}.
 The second is a real non-linear single  elliptic
equation for the length function of the metric  which comes from
 the Gauss curvature equation.

These equations can be approached
via a number of techniques in partial differential equations.  In this
short note, we improve upon  results obtained by assuming the Riemann surface
structure and postulating a fixed quadratic differential  representing the
$(2,0)$ part of a second fundamental form as solving the Codazzi
equations.   This leaves the problem of solving the elliptic Gauss equation
for the length function of the metric.  By analyzing the Gauss and
Codazzi equations together, we are able to reformulate the equations in
a form which completely identifies all local solutions in the case of negative
curvature. We prove that  the moduli space of solutions to the 
Gauss-Codazzi equations for a Riemann surface of genus greater than one 
immersed
with mean curvature  constant and less than $1$, in a not necessarily
complete  3-manifold of
constant curvature $-1$ is parameterized by cohomology classes of $(0,2)$
differentials.

The result is similar and proved in the same fashion as the results in gauge
theory in a paper of the first author \cite{G}. In fact, the details of how the
computations change with the change in base-point $g$ is comlicated, but it
is familiar to geometers from the variational formulation of the Yamabe 
problem and will not be repeated.   An abstract proof could be
constructed along the lines of the convexity theory used in the gauge
theory literature to describe bundle extensions by authors
Bradlow-Garcia-Prada \cite{B-G} and 
Daskalopolous-Uhlenbeck-Wentworth \cite{D-U-W}. 
Also, as presented here, the construction is not as natural as it would be
if viewed from the point of view of these authors. However,  the convexity
theory fails in cases of positive curvature, whereas the reformulation of
the Gauss-Codazzi equations as the Euler-Lagrange equations of a single
variational problem   has potential to contribute to the more
interesting cases of zero or even positive curvature.

It is entertaining 
to note that the second author came across these equations more than
twenty years ago in looking at the possibility that minimal surfaces could
be used to parameterize  quasi-Fuchsian hyperbolic 3-manifolds \cite{U}.
Although there is a fairly good existence theory, it seems as if
uniqueness is unlikely \cite{V}.

\section{The Variational formulation}

We  assume in this note that  $X$ is a surface immersed with
mean curvature $c$ in a 3-dimensional  manifold $N$ of constant sectional
curvature $k$. We do not assume this 3-manifold is complete.
The induced metric and second fundamental form of 
$X \subset  N$ satisfy
a system of Gauss-Codazzi equations which can be analyzed in a series of
steps as follows.

1.   
The induced metric on $X$ is of the form $h = h_{z\bar z} dz\,d\bar z$ 
where $(z,\bar z)$ are local 
complex coordinates in a complex structure $X_\sigma$  on $X$.

2.  
The second fundamental form $\gamma$      has the structure
$$\gamma = \alpha_{zz}(dz)^2 + ch_{z\bar z} dz\,d\bar z + \alpha_{\bar z\bar z}
(d\bar z)^2\ .$$
Here $c$ is the (constant) mean curvature. $\gamma$ is symmetric with respect
to the metric $h$, which implies the relationship 
$$\alpha_{zz} = \bar\alpha_{\bar z\bar z}\ .$$
         
3.  
The (2,0) part of $\gamma$ is a holomorphic quadratic differential 
$\alpha$ on $X$.

4.  
The induced (scalar) curvature $K$ of $X$ satisfies the Gauss equation.
$$K = k+c^{2} - 2|\alpha|^{2} .$$
Here the norm of $\alpha$   is assumed to be taken with respect to 
the metric $h$.
By  rescaling, we may assume that $k = (-1,0,1)$.  Our main  results
pertain to the case    $\lambda = k+c^2 <0$.

The solution to steps 1--3 is completely understood, so an automatic 
procedure would be  to fix the Riemann surface, the curvatures and the (2,0)
part alpha  of the second fundamental form and attempt to solve the Gauss 
equation for the metric.  There is a variational formulation of this
problem. This last step results in an equation which can be solved for
small $\alpha$ by the techniques of Kazden and Warner \cite{K-W}.
However, a slightly different variational problem  arises when
we solve for the second fundamental form and the metric in one step.

For convenience, we combine in a suggestive notation  $\lambda = k+c^2$. 

Fix the Riemann surface.  We will choose as a base point the constant
curvature metric on $X$ and a holomorphic quadratic differential. Let
$\beta_O$ be the $(0,2)$ form on $X$ which is dual to the chosen holomorphic
differential in the constant curvature metric.
Note that the identification between the $(2,0)$ and $(0,2)$ form depends
on the metric.  Now we compute with respect to an arbitrary element
  $\beta = \beta_0 + \bar\partial f_0$ and fixed  metric $g$ whose
conformal class determines the Riemann surface. Let $K(g)$ denote the
Gauss curvature of $g$.  Calculations are always easier at the base
metric, so we will want to have the freedom of changing it.

We call our
functional $\D$ for Donaldson, as it is in  reality a form  of
the Donaldson functional which appears  in the construction of
Hermitian-Einstein metrics in holomorphic vector bundles.  Usually 
there is no   explicit formula for this functional.  Due to the fact
that we are in line   bundles and are looking at an abelian gauge theory,
the functional is for us explicit. Note that it is well-defined up to
an overall constant, which we have fixed by assuming the functional $\D$
is $0$ on the constant curvature metric and the dual to the chosen
holomorphic quadratic differential. Note that the holomorphic quadratic
differential constructed from  the variational principle is not the one
that we started with. 

\begin{prop}\label{prop:gauss-codazzi}
The metric and holomorphic quadratic differential pair
$$(h,\alpha) = (e^{2u}g ,e^{2u} (\beta^*+\partial f^*))$$
solve the Gauss-Codazzi system (3--4)
if and only if  $u:x\to c$, $f:x\to  T^{1,0}x$ are critical maps for the
functional:
$$\D (u,f) = \iint \left[|\partial u|^2 + K(g)u + e^{2u} 
(-\lambda/2 + |\beta + \bar\partial f|^2)\right]\,d\mu + C(g)\ .$$
\end{prop}

\begin{proof}
This is a calculation. All the metrics, the covariant derivative     
$\partial = \bar\partial^*$ and the density are 
computed using the base 
metric $g$, although
 we suppress   this in the statement of the theorem. We have included
 the constant $C(g,\beta)$,   since  we wish to make this
 computation $(g,\beta)$ independent. The value of $C(g)$ is determined
by finding the value of the Donaldson functional at $g$ using our
original choice of constant curvature metric and holomorphic quadratic
differential.  

Now the equation which arises  from varying $f$ is a simple linear equation
(in $f$, not $u$)
$$\bar\partial \Big( e^{2u} (\beta^* +\partial f^*)\Big) =0\ .$$
This yields the holomorphic quadratic differential
$$\alpha = e^{2u} (\beta^* + \partial f^*)\ .$$
The equation obtained by varying $u$ is the equation
$$0 = - \partial^*\partial u + K(g) + e^{2u}( - \lambda + 2|\beta +
\bar\partial f|^2)\ .$$
Recall that in the new metric $h = e^{2u}g$  the curvature $K(h) =
e^{-2u} (K(g) - 2\partial\bar\partial u)$, 
so that the equation   which  arises when varying $u$ 
is indeed the desired equation for mean curvature. It can be written as 
$$K(h) = \lambda - 2e^{-4u}|\alpha|_g^{2}\ .$$ Note that the norm of
$\alpha$ is correctly computed in the new metric $h$.

\end{proof}

\begin{cor}
Fix $\Vol(X) = \iint  e^{2u}\, d\mu =T$.  Then the critical points of
 $$\hat \D(u,f) = \iint  |\partial u|^2 + k(g) u +  e^{2u} 
 |\beta + \bar\partial f|^2\,d\mu + C(g,\beta)$$
with respect to the constraint provide solutions of the Gauss-Codazzi
equations (3--4) with an unknown Lagrange multiplier   $\lambda$. 
\end{cor}

The actual utilization of this minimization principle seems quite delicate.
One would need  to employ the Moser inequality carefully.  Moreover, in many
cases one would be interested in saddle points rather than minima. However,
this variational principle which fixes the cohomology class of   
$\beta$ has definite 
advantages over the one which fixes the holomorphic differential $\alpha$.

\section{Solutions for $\lambda < 0$}

The main result of this short note  is the following theorem:

\begin{thm}\label{thm:main}
If $\lambda = k+c^2 <0$, there exists a unique solution  to the
 Gauss   Codazzi equations (3--4) for a fixed Riemann surface 
 and a (0,2) cohomology  class $[\beta] = \{\beta + \bar\partial f\}$. 
 \end{thm}

\begin{proof}  
We will show that for every solution to the Gauss Codazzi equations
with $\lambda  = k+c^2 < 0$, the Hessian  of $\D$ is positive definite.
For small $[\beta]$, there will be a solution near a constant negative 
curvature metric of curvature $\lambda$ on $X$ 
which can be found using an implicit
function theorem.  Openness follows from the invertibility of the Hessian.
Closedness is a rather easy calculation which we leave to the reader,
It is  similar to proofs in the literature \cite{K-W},\cite{U}.

Uniqueness follows since every solution can be connected to one with
small cohomology class  representative  $[\beta]$.This leaves the important
 step of showing the positive definiteness of the
Hessian to finish the proof.
\end{proof}

\begin{thm}\label{th:hessian}
If $(u,f)$ is a critical point of $\D$, for a fixed Riemann surface
and cohomology class $[\beta]$ with $\lambda < 0$, 
then the Hessian is positive definite.
\end{thm}

\begin{proof}
We might as well make the calculation with the solution metric and
holomorphic quadratic differential $\beta^*=\alpha$ as base-point. 
(Here is where our comments about basepoint pay off. 
These calculations  in gauge theory are where the idea for the proof comes
from).  

Since $(0,0)$ is a critical point, we have 
$$K(g) = \lambda - 2|\beta|^2$$ 
and 
$$\partial \beta=0\ .$$
The Hessian is easy to compute. 
$$H(v,f) = 2\iint  \left[ |\partial v|^2 - v^2 K(g) 
+ 4v \Re \langle  \beta, \partial f^*\rangle 
+ \langle \bar \partial  f,\partial f^*\rangle\right]\,d\mu\ .$$
We note that $\partial\beta =0$, so we may replace 
$$2\iint  v\Re \langle \beta,\partial f^*\rangle\,d\mu 
= -2 \iint \Re \langle \partial v\otimes f^*,\beta\rangle\,
d\mu\ .$$
Also 
$$\frac12 \iint \langle \bar\partial f,\partial f^*\rangle\,d\mu 
= \frac12\iint \langle\bar\partial f,\partial f^*\rangle 
- K(g) \langle f,f^*\rangle\,d\mu\ .$$

We have $K(g) = \lambda - 2|\beta|^2$ where $\lambda <0$. 
We rewrite the Hessian as 
\begin{equation*}
\begin{split} 
H(v,f) & = 2\iint \bigg[ |\partial v|^2 +  v^2 (-\lambda + 2|\beta|^2) 
+ 2v\Re \langle \beta,\partial f^*\rangle - 2\Re \langle \partial u\otimes 
f^*,\beta\rangle\\
&\qquad  + \frac12 \left( |\bar\partial f|^2 +\frac12 |\partial f|^2 
+ |f|^2 (-\lambda + 2|\beta|^2)\right) \bigg]\,d\mu\cr
& > 2\iint \bigg[ |\partial v|^2 + 2|v|^2 |\beta|^2 
+ \frac12 |\bar\partial f|^2 + |f|^2 |\beta|^2
- 2|v|\, |\beta|\, |\bar\partial f| 
- 2|\partial v|\, |f|\, |\beta|\bigg] \,d\mu \cr
&\ge 0\ .
\end{split}
\end{equation*}
Note that we may allow $\lambda =0$ in all but extremely degenerate cases.
\end{proof}

\end{document}